\documentclass[11pt,oneside]{amsart}
\usepackage{amscd}
\usepackage{amssymb}
\usepackage{amsfonts}
\usepackage[centertags]{amsmath}
\usepackage{latexsym}
\usepackage{verbatim}
\usepackage{amsthm}
\numberwithin{equation}{section}

\theoremstyle{plain}
\newtheorem{theorem}{Theorem}[section]

\newtheorem{lemma}[theorem]{Lemma}
\newtheorem{prop}[theorem]{Proposition}
\newtheorem{cor}[theorem]{Corollary}
\newtheorem{rem}[theorem]{Remark}
\newtheorem{defn}[theorem]{Definition}
\theoremstyle{definition}
\newtheorem{ex}{Example}

\renewcommand{\a}{\alpha}

\newcommand\C{{\mathbb C}}
\newcommand\R{{\mathbb R}}

\newcommand{\ov}[1]{\overline{ #1}}

\newcommand{\N}{\nabla}

\newcommand{\g}{\mathfrak{g}}

\begin{document}

\title[Tamed Symplectic forms and SKT metrics]{Tamed Symplectic forms and SKT metrics}
\author[N. Enrietti, A. Fino, L. Vezzoni]{Nicola Enrietti, Anna Fino and Luigi Vezzoni}
\address{Dipartimento di Matematica G. Peano \\ Universit\`a di Torino\\
Via Carlo Alberto 10\\
10123 Torino\\ Italy} \email{nicola.enrietti@unito.it}
 \email{annamaria.fino@unito.it}
 \email{luigi.vezzoni@unito.it}
\subjclass{32J27, 53C55, 53C30, 53D05}

\begin{abstract} Symplectic  forms  taming complex structures on compact manifolds
are strictly related to Hermitian metrics having the fundamental form $\partial \overline \partial $-closed, i.e.  to  strong K\"ahler with torsion  (${\rm SKT}$) metrics. It is still an open problem to exhibit a compact example of a complex manifold  having a tamed symplectic structure but non-admitting K\"ahler structures.
We show some negative results for the existence of   symplectic forms  taming  complex structures on  compact quotients of Lie groups by discrete subgroups. In particular, we prove that if $M$ is   a nilmanifold  (not a torus)  endowed with an invariant complex structure $J$, then   $(M, J)$  does not admit any symplectic
form taming $J$. Moreover, we show  that if a nilmanifold $M$ endowed with  an invariant complex structure $J$   admits an ${\rm SKT}$ metric, then $M$
is  at most  $2$-step. As a consequence we classify $8$-dimensional  nilmanifolds  endowed with an invariant complex structure  admitting an SKT metric.
\end{abstract}

\maketitle

\section{Introduction}
Let $(M,\Omega)$ be a compact $2n$-dimensional symplectic manifold. An almost complex structure $J$ on  $M$ is said to be \emph{tamed} by $\Omega$ if
$$
\Omega(X,JX)>0
$$
for any non-zero vector field $X$ on $M$.  When $J$ is a complex structure (i.e. $J$ is integrable) and $\Omega$ is tamed by $J$, the pair $(\Omega, J)$ has been called a {\em Hermitian-symplectic structure} in \cite{ST}.
Although any symplectic structure always admits tamed almost complex structures,
it is still an open problem to find an example of a compact complex manifold admitting a
Hermitian-symplectic structure, but no K\"ahler structures. From \cite{ST,LiZhang} there exist no examples in dimension $4$.
Moreover, the study of    tamed  symplectic structures in dimension $4$ is related to a more general conjecture of Donaldson (see for instance  \cite{donaldson,weinkove,LiZhang}).

A natural class  where to search examples of Hermitian-symplectic manifolds is provided by compact quotients of
Lie groups by discrete subgroups, since this set typically contains examples of manifolds admitting
both complex structures  and symplectic structures, but no K\"ahler structures.  Indeed,  it is well known that a {\em nilmanifold}, i.e. a compact
quotient of a nilpotent Lie group by a discrete subgroup, cannot
admit any K\" ahler metric unless it is a torus (see for
instance \cite{BG,Has}). Moreover, in the case of {\em solvmanifolds}, i.e.  compact
quotients of solvable  Lie groups by  discrete subgroups,
Hasegawa showed in \cite{Ha,Ha2} that a solvmanifold has  a
K\"ahler structure if and only if it is covered by a finite quotient
of a complex torus, which has the structure of a complex torus
bundle over a complex torus.

In this paper we show some negative results for the existence of Hermitian-symplectic structures on compact quotients of Lie groups by discrete subgroups.
First of all we show that Hermitian-symplectic structures are strictly related to a  \lq \lq special\rq\rq type of Hermitian metric, the SKT one. We recall (see \cite {GHR}) that a $J$-Hermitian metric $g$ on a complex manifold $(M,J)$ is called SKT (\emph{Strong K\"ahler with Torsion}) or \emph{pluriclosed}
 if the fundamental form $\omega(\cdot,\cdot)=g(J\cdot,\cdot)$ satisfies
$$
\partial\ov \partial \omega=0\,.
$$
For complex surfaces a  Hermitian metric satisfying the SKT condition is  {\emph{standard}} in the terminology of Gauduchon \cite{Gau2} and one standard metric  can be found in the conformal
class of any given Hermitian metric on a compact manifold. However,  the theory is different in higher dimensions.
SKT metrics have a central role in type II string theory, in $2$-dimensional supersymmetric
$\sigma$-models (see \cite{GHR,strominger})  and they have also relations with generalized
K\"ahler geometry (see for instance
\cite{GHR,Gu,Hi2,AG,CG,FT0}).

SKT metrics have been recently studied by many
authors. For instance, new simply-connected  compact strong KT examples have been recently
constructed by Swann in \cite{Sw} via the twist construction. Moreover, in  \cite{FT1} it has been shown that the blow-up  of an SKT  manifold $M$ at a point or along a compact submanifold  admits an SKT
metric.

For real Lie groups admitting SKT metrics there are only some classification results in  dimension $4$ and $6$. More precisely,
$6$-dimensional SKT nilpotent Lie groups have been classified  in \cite{FPS}  and  a classification of SKT
solvable Lie groups of dimension $4$ has been recently obtained in \cite{MS}.
\medskip

We outline our paper. In Section 2 we show that  the existence of a Hermitian-symplectic structure  on a complex manifold  $(M, J)$ is equivalent to the existence of a $J$-compatible SKT metric $g$ whose fundamental form $\omega$  satisfies
$\partial\omega=\ov\partial\beta$ for some $\partial$-closed  $(2,0)$-form $\beta$ (see Proposition \ref{linkwithSKT}).
This result allows us to write down a natural obstruction to the existence of a Hermitian-symplectic structure on a  compact complex
manifold (see Lemma \ref{fond}).

Let $G/\Gamma$   a compact quotient of a simply-connected Lie group $G$  by a discrete subgroup $\Gamma$.   By  an invariant  complex   (resp. symplectic)  structure on  $G/\Gamma$ we will mean  one induced by  a complex   (resp. symplectic)  structure  on  the Lie algebra $\mathfrak g$  of $G$.
In Section 3 we start  by showing that the existence of a
symplectic form  taming an invariant complex structure  $J$ on $G/\Gamma$   implies the existence of an invariant one (see Proposition \ref{noninvHS}). This can be used to prove the following

\begin{theorem}\label{main1}
Let $(M=G/\Gamma,J)$ be a compact quotient of a Lie group by a discrete subgroup  endowed with an  invariant complex structure $J$.
Assume that the center $\xi$ of the Lie algebra $\g$ associated to $G$ satisfies
\begin{equation}\label{Jxi}
J\xi\cap [\g,\g]\neq \{0 \}\,,
\end{equation}
then $(M,J)$ does not admit  any symplectic form taming $J$.
\end{theorem}

In the case of a a simply-connected
nilpotent Lie group $G$ if the structure equations of its Lie algebra $\mathfrak g$ are rational,
then  there exists a discrete subgroup $\Gamma$  of $G$ for which the quotient $G/\Gamma$  is compact \cite{Malcev}. A general classification of real nilpotent Lie algebras exists only for dimension less or  equal than  $7$ and $6$ is the highest dimension in which there do not exist continuous families (\cite{Magnin, GK}).
It turns out that in dimension $6$  by \cite{FPS}  any $6$-dimensional non-abelian
nilpotent Lie algebra admitting an SKT metric is $2$-step nilpotent. Moreover, by \cite{Sw} any  $6$-dimensional SKT nilmanifold  is the   twist of the K\"ahler product of the two torus ${\mathbb T}^4$ and ${\mathbb T}^2$ by using  two integral $2$-forms  supported on ${\mathbb T}^4$.
In this paper we extend the  previous result     to any dimension.
\begin{theorem}\label{main2}
Let $(M=G/\Gamma,J)$ be a  nilmanifold  $($not a torus$)$ endowed with an invariant complex structure $J$. Assume that there exists a $J$-Hermitian {\rm SKT} metric $g$ on $M$. Then ${\mathfrak g}$ is $2$-step nilpotent and $M$  is a  total space of a  principal holomorphic torus bundle over a torus.
\end{theorem}

Moreover, we show that an invariant  SKT metric  compatible with an invariant complex structure on a nilmanifold is always standard (Corollary \ref{standard}).

\smallskip

Theorem \ref{main1} and Theorem \ref{main2} will be used to prove the following
\begin{theorem}\label{main3}
Let $(M=G/\Gamma,J)$ be a nilmanifold  $($not a torus$)$ with an invariant complex structure. Then $(M,J)$ does not admit any symplectic form taming $J$.
\end{theorem}
It can be observed that there exist examples of nilpotent Lie algebras not satisfying condition \eqref{Jxi}
(see Example \ref{ex10}). Hence Theorem \ref{main3} cannot be directly deduced  from Theorem \ref{main1}.

 In dimension $8$ some examples of nilpotent Lie algebras endowed with an SKT metric have been found in \cite{Rasmus}, \cite{FT2} and \cite{RT}, but there is no general classification result.
In the last part of the paper we describe  the $8$-dimensional  nilmanifolds  endowed with an invariant complex structure  admitting an SKT metric.
In dimension $6$  by \cite{FPS}  the existence of a
strong KT structure on a nilpotent Lie algebra $\mathfrak g$  depends only on the complex structure of  $\mathfrak g$.
 We show   that in dimension $8$ the previous property  is  no longer true.

A good quaternionic analog of K\"ahler geometry is given by {\em hyper-K\" ahler  with
torsion} (shortly HKT) geometry.  This geometry
was introduced by Howe and Papadopoulos  \cite{HP} and later studied for instance
in \cite{GP}.  Nilmanifolds  also provides examples of compact HKT manifolds (see \cite{GP,
DF}).  In the last section we investigate which $8$-dimensional SKT nilmanifolds admit also HKT structures.

\medskip

\noindent \textbf{Acknowledgements.} We would like to thank  Paul Gauduchon, Thomas Madsen, Andrew Swann and Adriano Tomassini for  useful conversations  and comments on the paper. We also would like to thank the referee for valuable
remarks  which  improved the contents of the paper. This work was supported by the Projects MIUR ``Riemannian Metrics and Differentiable Manifolds'',
``Geometric Properties of Real and Complex Manifolds'' and by GNSAGA
of INdAM.

\smallskip

\section{Link with SKT metrics}

The study of SKT metrics  is strictly related to the study of the geometry of the Bismut connection. Indeed, any Hermitian structure $(J,g)$ admits a unique connection $\N^B$ preserving $g$ and $J$ and such that the tensor
$$
c (X, Y, Z) = g(X, T^B(Y , Z ))
$$
is totally skew-symmetric, where $T^B$ denote the torsion of $\N^B$ (see \cite{Gau}).
This connection was used by Bismut in \cite{Bismut}  to prove a local index
formula for the Dolbeault operator for non-K\"ahler manifolds.
The torsion $3$-form $c$ is related to the fundamental form $\omega$  of $g$ by
$$
c(X, Y, Z)  =  - d \omega (JX, JY, JZ)
$$
and it is well known that
$
 \partial\ov \partial \omega=0$ is equivalent to  the condition $$ dc=0\,.
$$

By  \cite{StreetsTian} it turns out that the Hermitian-symplectic structures are related  to  static solutions of a   new metric flow on complex manifolds called
\emph{Hermitian curvature flow}. Indeed,
 Streets and Tian constructed  a flow using the Ricci tensor associated to the Chern connection
instead of the Levi-Civita connection. In this way they obtained an elliptic flow and proved some results on
short-time existence of solutions for this flow and on stability of K\"ahler-Einstein metrics. A modified Hermitian
curvature flow was used in \cite{ST} to study the  evolution of  SKT metrics, showing  that the existence of  some particular type of   static SKT metrics  implies the existence of a Hermitian-symplectic structure on the complex manifold. Static SKT metrics on Lie groups have been also recently studied in \cite{Enrietti}.

The aim of this section is to point out a natural link between Hermitian-symplectic structures  (also called \emph{holomorphic-tamed} by de Bartolomeis and Tomassini \cite{Tomas})  and SKT metrics. We have the following
\begin{prop}  \label{linkwithSKT} Let $(M, J)$ be a complex manifold.
Giving a Hermitian-symplectic structure  $\Omega$ on  $(M, J)$ is equivalent to assigning an {\rm SKT} metric $g$ such that the associated fundamental form $\omega$ satisfies
\begin{equation}\label{omega}
\partial \omega=  \overline{\partial}\beta
\end{equation}
for some $\partial$-closed  $(2,0)$-form $\beta$.
\end{prop}
\begin{proof}
Let $\Omega$ be a Hermitian-symplectic structure on $(M, J)$, then:
$$d\Omega=0,  \quad \Omega^{1,1}>0,$$
where $\Omega^{1,1}$  denotes the $(1,1)$-component of $\Omega$.
 We can write the real $2$-form $\Omega$ as
$$
\Omega=\omega-\beta-\overline{\beta}\,.
$$
where $\beta=-\Omega^{2,0}$ is the opposite of the $(2,0)$-component of $\Omega$ and  $\omega=\Omega^{1,1}$ is the $(1,1)$-component. Then
$$
d\Omega=0\Longleftrightarrow
\begin{cases}
\partial\omega=\overline{\partial}\beta\\
\partial\beta=0\,
\end{cases}
$$
and $\omega$ is the fundamental  form associated  to an SKT metric.

Conversely if $g$ is an SKT metric whose  fundamental form $\omega$ satisfies \eqref{omega} for some  $\partial$-closed $(2,0)$-form  $\beta$, then
$\Omega=\omega-\beta-\overline{\beta}$ defines a Hermitian-symplectic structure on $(M,J)$.
\end{proof}
In order to write down a natural obstruction to the existence of a Hermitian-symplectic structure, we recall some basic facts about Hermitian Geometry.
Let $(M,J,g)$ be a $2n$-dimensional compact Hermitian manifold and let  $\omega$  be the associated fundamental form.
The complex structure $J$  is extended to act on $r$-forms by
$$
J\alpha(X_1,\dots,X_r)=(-1)^r \alpha(JX_1,\dots,JX_r)\,.
$$
With respect to this extension $J$ commutes with the duality induced by $g$. On $M$  the Hodge  star operator
$$
*\colon \Lambda^{r,s}(M)\to \Lambda^{n-s,n-r}(M)
$$
is defined  by the relation
$$
\alpha\wedge*\ov{\beta}=g(\alpha,\ov{\beta})\,\frac{\omega^n}{n!}\,.
$$
Moreover, $g$ induces the $L^2$ product
$$
(\alpha,\beta)=\int_M \alpha\wedge *\ov{\beta}\,\frac{\omega^n}{n!}
$$
and the differential operators
$$
\partial^{*}\colon \Lambda^{r,s}(M)\to \Lambda^{r-1,s}(M)\,,\quad \ov{\partial}^{*}\colon \Lambda^{r,s}(M)\to \Lambda^{r,s-1}(M)
$$
defined as
$$
\partial^{*}=-*\ov{\partial}\,*\,, \quad \ov{\partial}^{*}=-*\partial*\,.
$$
It is well known that
$$
(\partial^*\alpha,\beta) =(\alpha,\partial\beta) \,,\quad (\ov{\partial}^*\alpha,\beta) =(\alpha,\ov{\partial}\beta) \,.
$$

\begin{lemma}
\label{fond}
Let $(M,\Omega,J)$ be a  compact Hermitian-symplectic manifold and let $\eta$ be a $(2,1)$-form. Fix an arbitrary Hermitian metric $g$ on  $(M,J)$ with associated Hodge star operator $*$ .
Then
$$
(\partial^*\eta,\Omega^{1,1})\neq 0\Longrightarrow \overline{\partial}^*\eta\neq 0\,.
$$
\end{lemma}
\begin{proof}
We have
$$
(\partial^*\eta,\Omega^{1,1})=(\eta,\partial \Omega^{1,1}).
$$
By Proposition \ref{linkwithSKT} $$\partial \Omega^{1,1} = \overline \partial \beta,$$
for some $\partial$-closed $(2,0)$-form. Therefore
$$
(\partial^*\eta,\Omega^{1,1}) =(\eta,\overline{\partial}\beta)=(\overline{\partial}^*\eta,\beta)
$$
which implies the statement.
\end{proof}

\section{Main results}
In this section we prove the results stated in the introduction.

We recall that an  almost complex structure $J$  on a real Lie algebra $\mathfrak g$  is
said to be \emph{integrable} if the Nijenhuis condition
$$
[X, Y] - [JX, JY]  + J[JX, Y] + J[X, JY] = 0
$$
holds for all $X, Y \in \mathfrak g$.  By a complex structure on a Lie algebra we will always mean an
integrable one. We can give the following
\begin{defn}
Let $({\mathfrak g},J)$ be a real Lie algebra endowed with a complex structure $J$. A {\rm Hermitian-symplectic}  structure on $({\mathfrak g},J)$ is a  symplectic  form
 $\Omega$ which tames $J$, that is  a real $2$-form $\Omega$ such that
$$
\Omega(X,JX)>0,
$$
for  any non-zero vector  $X \in \mathfrak g$.
\end{defn}
In order to prove Theorem \ref{main1}, we consider the following
\begin{lemma} \label{MainLieAlg} Let ${\mathfrak g}$ be a real  Lie algebra endowed with a  complex structure $J$  such that $J \xi \cap [\mathfrak g, \mathfrak g] \neq \{ 0 \}$,
where  $\xi$ denotes the center of $\mathfrak g$.
Then $({\mathfrak g}, J)$ cannot admit any Hermitian-symplectic structure.
\end{lemma}

\begin{proof} Suppose that $({\mathfrak g}, J)$ admits a Hermitian-symplectic structure $\Omega$, then by definition $\Omega$  satisfies the  conditions
$$
d\Omega (X, Y, Z) = - \Omega([X, Y], Z) - \Omega([Y, Z], X])
 - \Omega ([Z, X], Y)=0,
$$
for every
$X, Y,  Z \in {\mathfrak g}$ and
$$
\Omega  (X, J X) > 0,  \quad  \forall \, X \neq 0.
$$
Therefore, in particular one has:
$$
d\Omega (X, Y, Z) =  - \Omega([Y, Z], X) =0, \quad \forall X \in \xi, \forall\, Y, Z \in \mathfrak g,
$$
or, equivalently, that
$$
\Omega( X, W) =0,
$$
for any $X \in \xi$ and $W \in [\mathfrak g, \mathfrak g].$   Assume that exists  a non-zero $W \in [\mathfrak g, \mathfrak g] \cap J \xi$; then $JW \in \xi$  and one has $\Omega (W, JW) =0$ which is not possible since  $J$ is  tamed by $\Omega$.

 \end{proof}

If $G$ is a real Lie group with Lie algebra $\mathfrak g$, then assigning  a left-invariant
almost complex structure on $G$ is equivalent to  choosing an almost complex
structure $J$ on $\mathfrak g$. Such a $J$  is integrable if and only if it is integrable as an almost
complex structure on $G$.  Therefore, a complex structure on $\mathfrak g$ induces a complex structure on $G$ by the
Newlander-Nirenberg theorem  and $G$ becomes a complex
manifold. The elements of $G$ act holomorphically on $G$  by multiplication on the left but $G$  is not a complex Lie group in general.

Let now $M = G/\Gamma$ be a compact quotient of a simply-connected Lie group $G$ by a uniform discrete subgroup $\Gamma$.
By  an {\em invariant} complex structure  (resp.  an {\em invariant} Hermitian-symplectic structure)  on $G/\Gamma$  we  will mean  a complex structure (resp. Hermitian-symplectic structure)
 induced by one on the Lie algebra ${\mathfrak g}$.

We can prove the following:

\begin{lemma} \label{noninvHS}
If $(M = G/\Gamma, J)$ with $J$ an invariant  complex structure  admits a Hermitian-symplectic
structure $\Omega$, then it admits an invariant   Hermitian symplectic structure $\tilde \Omega$.
\end{lemma}

\begin{proof} Suppose that there exists a non-invariant Hermitian-symplectic structure  $\Omega$, then by using the property that $G$ has a bi-invariant volume form $d \mu$ (see \cite{Mi}) and by applying the symmetrization process of \cite{FG} we can construct a new symplectic invariant form  $\tilde \Omega$, defined by
$$
\tilde{\Omega}(X,Y)=\int_{m \in M}  \Omega_m(X_m,Y_m)\, d \mu,
$$
for any left-invariant vector field $X, Y$. The $2$-form $\tilde \Omega$  tames the complex structure $J$.
\end{proof}

Now are ready to prove Theorem \ref {main1}.

\begin{proof}[Proof of Theorem $\ref{main1}$]
Let $(M=G/\Gamma,J)$ be a compact quotient of a Lie group with a left-invariant complex structure. Assume that there exists a $J$-compatible
Hermitian-symplectic structure on $M$. Then using Lemma \ref{noninvHS}, we can construct a left-invariant Hermitian-symplectic structure $\Omega$
on $(G,J)$.  In view of Lemma \ref{MainLieAlg}, we have $J\xi\cap [\g,\g]= \{ 0 \}$, as required.
\end{proof}

Consider a Lie algebra $(\g,J,g)$ with a complex structure and a $J$-compatible inner product.
Let $\N^B$ be the Bismut connection associated to $(J,g)$.
Then  the inner product $g$ is ${\rm SKT}$  if and only if the torsion $3$-form $c(X,Y,Z)=g(X,T^{B}(Y,Z))$ of the Bismut connection $\nabla^B$  is closed. In view of \cite{DF},  the Bismut connection
$\N^B$ and
the torsion  $3$-form $c$ can be written in terms of Lie brackets as
\begin{equation}\label{B}
\begin{split}
(\nabla^B_X Y, Z) =\, \frac{1}{2} \{&\, g ([X, Y] -[JX, JY], Z) - g ([Y, Z] +[JY, JZ], X) \\
&+ g ([Z, X] -[JZ, JX], Y)\}
\end{split}
\end{equation}
\begin{equation}\label{c}
c(X, Y, Z) = - g([JX, JY], Z)  - g([JY, JZ], X) -  g([JZ, JX], Y).
\end{equation}
for any $X,X,Z\in \g$. We have the following
\begin{lemma} \label{difftorsion} Let $(\g,J,g)$ be a Lie algebra with a complex structure and a $J$-compatible inner product.
Let $\xi$ be the center of $\g$ and let $X\in\xi$. Then for any $Y\in \g$
\begin{equation}\label{dc}
dc (X, Y, JX, JY) = 2\Big( \|[Y, JX]\|^2-   g([[JX, Y],JX], Y)  -  g([[Y, JY], JX], X)\Big)\,.
\end{equation}
\end{lemma}
\begin{proof}
Let $X\in\xi$. Then since $J$ is integrable and $X$ belongs to the center, we have
\begin{equation} \label{intdeemcenter}
[JX, JY] = J [JX, Y].
\end{equation}
Using \eqref{B} and \eqref{c} we have
\[
\begin{split}
dc (X, Y,&\, JX, JY)= \\
=& - c([Y, JX], X, JY) + c([Y, JY], X, JX) - c([JX, JY],X, Y)\\[3 pt]
=& - c([Y, JX], X, JY) + c([Y, JY], X, JX) - c(J [JX, Y],X, Y)\\[3 pt]
=&\ g ([J [Y, JX], JX], JY) + g ([Y, JX], [Y, JX]) - g([Y, J [Y, JX]], X)\\
& - g([J[Y, JY], JX], JX)- g([[JX, Y],JX], Y) \\
& + g ([JX, JY], [JX, JY]) - g ([JY, [JX, Y]], X).
\end{split}
\]
Now, by \eqref{intdeemcenter} we get $[J [Y, JX], JX] = J [[Y, JX], JX]$, so
$$
g([J[Y,JX],JX],JY)-g([[JX, Y], JX],Y)=-2g([[JX, Y],JX], Y).
$$
Moreover
\[
\begin{split}
g([Y, J [Y, JX]], X) +g([J[Y, JY], JX], JX) +g& ([JY, [JX, Y]], X) =\\[3pt]
  &\ =2 g([[Y, JY], JX], X)
\end{split}
\]
and therefore
\[
\begin{split}
dc (X,&\, Y, JX, JY) =\\[3pt]
=&\, 2 g ([Y, JX], [Y, JX])- 2  g([[JX, Y],JX], Y)  - 2 g([[Y, JY], JX], X)\,,
\end{split}
\]
as required.
\end{proof}

Now we consider the nilpotent case. We recall that a  real Lie  algebra $\mathfrak g$ is nilpotent  if there exists a  positive $s$ such that for  the descending central series $\{ {\mathfrak g}^k \}_{k \geq 0}$ defined by
$$
{\mathfrak g}^0 = {\mathfrak g},  \quad {\mathfrak g}^1= [{\mathfrak g}, {\mathfrak g}]\,,\quad\ldots\,,\quad
{\mathfrak g}^{k }  =[{\mathfrak g}^{k - 1}, {\mathfrak g}]\,,
$$
one has  ${\mathfrak g}^{s}= \{0 \}$ and ${\mathfrak g}^{s-1} \neq \{0 \}$. The Lie algebra ${\mathfrak g}$ is therefore called $s$-step nilpotent  and one has that ${\mathfrak g}^{s -1}$ is contained in the center $\xi$ of $\mathfrak g$.

The next proposition gives an obstruction to the existence of an SKT inner product on
a nilpotent Lie algebra with a complex structure.

\begin{prop} \label{centerjinv}
Let $({\mathfrak g}, J)$ be a nilpotent Lie algebra endowed with a complex structure $J$. If the center $\xi$ is not $J$-invariant, then $({\mathfrak g}, J)$ does not admit  any  {\rm SKT} metric.
\end{prop}

\begin{proof}
Let $g$ be an arbitrary $J$-compatible inner product on $\g$ and let  $\nabla^B$ be the associated Bismut connection. Assume that the center $\xi$ is not $J$-invariant, then
there exists $X \in \xi$ such that  $J X\notin\xi$, or equivalently there exists $X\in\xi$ such that the map
$ad_{JX}: {\mathfrak g} \longrightarrow {\mathfrak g}$ is not identically zero. But certainly the restriction $ad_{JX} \vert_{{\mathfrak g}^{s-1}}$ is zero because $\mathfrak g^{s-1} \subseteq \xi$. So there exists an integer $k$ such that $ad_{JX} \vert_{{\mathfrak g}^{k}}$ is not identically zero and $ad_{JX} \vert_{{\mathfrak g}^{k+1}}$ is zero. Then we can choose $Y\in \mathfrak g^k$ such that $[JX,Y]\neq 0$ and $[[Y,Z],JX]=0$ for each $Z\in\mathfrak g$; so,  using \eqref{dc} we have
\[
dc(X, Y, JX,JY)  =2\,\|[Y, JX]\|^2 \neq 0
\]
and $g$ cannot be {\rm SKT}.
\end{proof}

Let $J$ be a complex structure on a $2n$-dimensional nilpotent Lie algebra $\mathfrak g$. We introduce the ascending series $\{ {\mathfrak g}^J_l \}_{l \geq 0}$ defined inductively by
$$
\begin{cases}
{\mathfrak g}^J_0 = \{ 0 \} \\
{\mathfrak g}^J_l  = \{ X \in {\mathfrak g} \mid [J^k X, {\mathfrak g} ] \subseteq {\mathfrak g}^J_{l - 1}, \,  k = 1,2 \},\quad l>1.
\end{cases}
$$
We say that $J$ is \emph{nilpotent} if ${\mathfrak g}^J_h = {\mathfrak g}$ for some positive integer $h$.
Equivalently by \cite{CFGU},  a complex structure on a $2n$-dimensional  nilpotent Lie algebra $\mathfrak g$ is nilpotent if there is a basis
$\{  \alpha^1, \dots, \alpha^n \}$ of $(1,0)$-forms
satisfying  $d \alpha^1 =0$ and
$$
d \alpha^j \in \Lambda^2 \langle \alpha^1, \ldots, \alpha^{j - 1}, \overline \alpha^1, \ldots, \overline \alpha^{j - 1} \rangle,
$$
for $j = 2, \ldots, n$.

Let $({\mathfrak g},J)$ be a nilpotent Lie algebra with a left-invariant complex structure admitting an SKT metric. Since the center $\xi$ of $\g$ is $J$-invariant, $J$ induces a complex structure $\hat J$ on the quotient $\hat \g={\mathfrak g}/{\xi}$ by the relation
\begin{equation} \label{indcxstruc}
\hat J (X + \xi) = JX + \xi
\end{equation}
for any $X\in\g$.
Moreover, if $\hat J$ is nilpotent, then $J$ is also nilpotent. Indeed, we have
$$
{\mathfrak g}_1^J = \xi, \quad   {\mathfrak g}_2^J = {\rm span} \langle \hat {\mathfrak g}_1^J, \xi\rangle
$$
and more in general $ {\mathfrak g}_k^J = {\rm span} \langle \hat {\mathfrak g}_{k - 1}^J, \xi \rangle $, for any $k \geq 1$.
Therefore, if  $J$ is non-nilpotent, $\hat J$ has to be non-nilpotent.

 We have the following
\begin{prop}\label{indotta} Let $({\mathfrak g}, J)$ be a nilpotent Lie algebra endowed with an {\rm SKT} inner product $g$.
Then there exists an {\rm SKT} inner product  $\hat g$ on $(\hat \g={\mathfrak g}/{\xi},\hat J)$, where  $\hat J$ is defined by \eqref{indcxstruc}.
\end{prop}

\begin{proof} We have the decomposition
$$
{\mathfrak g} = \xi \oplus  \xi^{\perp},
$$
where $\xi^{\perp}$ denotes the orthogonal complement of $\xi$ with respect to the SKT metric $g$. So any $X\in\g$ splits in
$$
X=X^{\xi}+X^{\perp}\,.
$$
Then we can identify at the level of vector spaces ${\mathfrak g}/{\xi}$ with $\xi^{\perp}$ in the following way
$$
{\mathfrak g}/{\xi} \cong  \{ X + \xi \, \mid \, X \in \xi^{\perp} \}.
$$
We claim that the inner product $\hat g$ on $\hat \g$ defined by
$$
\hat g (X + \xi, Y + \xi) = g(X^{\perp}, Y^{\perp}), \quad X,Y\in\g
$$
is SKT. We note that $[X,Y]^\perp=[X^\perp,Y^\perp]$, and applying Proposition \ref{centerjinv} we obtain $(JX)^\perp=JX^\perp$. Then by a direct calculation the torsion $3$-form of the Bismut connection of $(\hat J, \hat g)$ satisfy
$$
\hat c(X + \xi , Y + \xi , Z + \xi)=c(X^{\perp}, Y^{\perp}, Z^{\perp})
$$
for any $X, Y, Z \in \g$.  Therefore the closure of $c$ implies the closure of $\hat c$ since $dc(X, Y, Z, W)=0$ for every $X, Y, Z, W \in {\xi}^{\perp}$.
\end{proof}

Now are ready to prove Theorem \ref{main2} and Theorem \ref{main3}.
\begin{proof}[Proof of Theorem $\ref{main2}$]
Let $(M= G/\Gamma,J)$ be a    nilmanifold with an invariant complex structure $J$. Assume that there exists a $J$-Hermitian SKT metric $g$ on $M$. We can suppose  that $g$ is left-invariant (see \cite{Ugarte,FG}).
Then $(J,g)$ can be regarded as an SKT structure on the Lie algebra $\g$ associated to $G$.
Since $J$ preserves the center $\xi$ of ${\mathfrak g}$ and $({\mathfrak g}/{\xi}, \hat J)$ admits an SKT metric, the center  of  ${\mathfrak g}/{\xi}$
$$
\hat \xi = \{ X^{\perp} + \xi \mid [X^{\perp}, \g] \in \xi  \},
$$
is $\hat J$-invariant. Then the ideal of $\g$
$$
\pi^{-1}\hat\xi= \{ X \in {\mathfrak g} \mid [X, {\mathfrak g}]  \subseteq \xi \}
$$
is $J$-invariant.

We prove that  ${\mathfrak g}$  is  $2$-step nilpotent by induction on the  dimension of $\g$. If  ${\rm dim}\,\g=6$, then  the claim follows from \cite{FPS}. Suppose now that the statement is true for $\dim \g< 2n$ and assume $\dim \g=2n$. Let $\xi$
be the center of $\g$ and let $\hat \g=\g/\xi$. In view of Proposition \ref{indotta}  we have that $(\hat \g,\hat J)$ admits an SKT inner product. Hence, since $\dim \hat \g < 2n$   by our assumption we get that   $\hat \g$ is at most $2$-step nilpotent.
This implies that ${\mathfrak g}$ is  at most $3$-step nilpotent. If  ${\mathfrak g}$ is $3$-step, then there exist $X, Y \in \mathfrak g$ such that
$$
[X, Y] \neq 0.
$$
with $Y =   [X_1, X_2] \in {\mathfrak g}^1 \subseteq \pi^{-1}\hat\xi$. Since $J(\pi^{-1}\hat\xi) = \pi^{-1}\hat\xi$, we have that $JY \in \pi^{-1}\hat\xi$. Now
\[
\begin{split}
dc(X, Y,&\, JX, JY)  = \\[3pt]
=&\,-c([X,Y], JX, JY) + c([X, JX], Y, JY) - c([X, JY], Y, JX)\\
&\,- c([Y, JX], X, JY) + c([Y, JY], X, JX) - c([JX, JY], X, Y)\\[3pt]
=&\, g ([X,Y], [X,Y]) - g([J[X,JX],JY], JY) + g([Y, J[X, JX]], Y)\\
&\, - 2 g([Y, JY], [X, JX])+ g([X, JY], [X, JY]) \\
&\, + g ([Y, JX], [Y, JX])+ g([JX, JY], [JX, JY]).
\end{split}
\]
As a consequence of the  Jacobi identity we get
\begin{equation} \label{Jacobig_2}
[[W_1, W_2], Z] =0,
\end{equation}
for any $W_1, W_2 \in {\mathfrak g}$ and $Z \in  \pi^{-1}\hat\xi$. Therefore  $[Y, JY]$ and $[Y, J[X, JX]]$ vanish. Now we can show also that $[J[X,JX],JY]$ is zero. Indeed,
by using  the integrability of $J$  and \eqref{Jacobig_2}   we get
$$
[J[X,JX],JY] = [[X,JX],Y] + J[J[X,JX], Y] + J [[X,JX], JY] = 0,
$$
and so
$$
0=dc(X, Y, JX, JY) =\| [X,Y] \|^2 + \| [X, JY] \|^2 +  \| [Y, JX] \|^2 + \| [JX, JY] \|^2\,,
$$
which is a contradiction with the assumption that $\mathfrak g$ is $3$-step nilpotent.

Since the center $Z(G)$ of $G$  is $J$-invariant and $\Gamma \cap Z(G)$ is a uniform discrete subgroup  of  $Z(G)$, we have that the suriective homomorphism $\pi: {\mathfrak g} \to {\mathfrak g}/{\xi}$ induces  a  holomorphic principal  torus bundle $\tilde \pi: G/\Gamma \to {\mathbb T}^{2p}$ over a $2p$-dimensional torus ${\mathbb T}^{2p}$ with  $2p = \dim \mathfrak g - \dim \xi$.
\end{proof}

\begin{rem}
{\rm  \begin{enumerate}

\item Since  by \cite{Rollenske}  on a $2$-step nilpotent Lie algebra every complex structure is nilpotent, Theorem \ref{main2} implies
that if  a nilmanifold endowed with a left-invariant complex structure $J$ admits  a $J$-compatible  SKT metric, then the complex structure has to be
nilpotent.

\item By \cite{Sw} any  $6$-dimensional SKT nilmanifold  is the   twist of the K\"ahler product of the two torus ${\mathbb T}^4$ and ${\mathbb T}^2$ by using  two integral $2$-forms  supported on ${\mathbb T}^4$. In a similar way  as  a consequence  of Theorem \ref{main2} we have then that any SKT nilmanifold  is the twist of a torus.
\end{enumerate}}
\end{rem}

\begin{cor} \label{standard} Let $(M=G/\Gamma,J)$ be a  $2n$-dimensional nilmanifold  endowed with an invariant complex structure $J$. If  there exists a $J$-compatible {\rm SKT} metric $g$ on $M$  and $n > 2$,  then $b_1(M) \geq 4$ and  every    invariant $J$-compatible    metric  on $M$ is standard.
\end{cor}

\begin{proof} By Theorem \ref{main2} we have the Lie algebra $\g$ of $G$ has to be $2$-step nilpotent. Since ${\mathfrak g}^1  \subseteq \xi$, the annihilator of  $\xi$ in ${\mathfrak g}^*$ is contained in the annihilator  of  ${\mathfrak g}^1 $, therefore its elements are all $d$-closed.  If $\dim \xi = 2n - 2p$ we can choose a basis of $(1,0)$-forms $\{\alpha^1, \ldots, \alpha^n \}$ such that
$$
\left \{ \begin{array}{l}
d \alpha^j = 0, \quad j = 1, \ldots, p,\\
d \alpha^l \in \Lambda^2 \langle \alpha^1, \ldots, \alpha^p, \overline \alpha^1, \ldots , \overline \alpha^p \rangle, \quad  l = p + 1, \ldots,n,
\end{array} \right.
$$
with $\{ \alpha^1, \ldots, \alpha^p \}$ in the annihilator of $\xi$ and $\{ \alpha^{p + 1}, \ldots, \alpha^n \}$ in the complexification of  the dual $\xi^*$ of $\xi$.
Using the previous  basis we can show that there exist at least two $d$-closed $(1,0)$-forms. Indeed,
 if   $p = 1$, then we have
$$
\left \{ \begin{array}{l}
d \alpha^1= 0, \\
d \alpha^l \in \Lambda^2 \langle  \alpha^1,  \overline \alpha^1 \rangle, \quad  l = 2, \ldots,n,
\end{array} \right.
$$
and therefore ${\mathfrak g}$ is isomorphic to the direct sum  ${\mathfrak h}_3^{\R} \oplus  \R^{2n - 3}$ of the $3$-dimensional real Heisenberg Lie algebra ${\mathfrak h}_3^{\R}$ with $\R^{2n - 3}$ and $b_1 ({\mathfrak g}) = 2n - 1 \geq4$ since $n>2$.

By using Nomizu's Theorem \cite{Nomizu} we have that the de Rham cohomology of $M$ is isomorphic to the Chevalley-Eilenberg cohomology of the Lie algebra ${\mathfrak g}$. Since there exist at least four real  $d$-closed $1$-forms we get that $b_1( M) = b_1 ({\mathfrak g}) \geq 4$.

We recall that by \cite{Gau2} a  $J$-compatible   metric $h$ on $(M, J)$  is standard if and only  if the Lee form  $\theta = J d^* \omega$  is co-closed, where $\omega$ denotes the fundamental form of $(J, h)$.  Since ${\mathfrak g}$ is nilpotent, it is also unimodular, so any $(2n -1)$-form is $d$-closed and in particular,  if $\theta$ is the Lee form of an invariant  $J$-compatible  metric $h$, then   $d(* \theta) =0.$
\end{proof}

\begin{proof}[Proof of Theorem $\ref{main3}$]
Let $(M,G/\Gamma,J)$ be a nilmanifold with an invariant complex structure and let $(\g,J)$ be its associated Lie algebra. Assume that there exists a Hermitian-symplectic structure $\Omega$ on $(M,J)$. Using Lemma \ref{noninvHS}, we may assume that $\Omega$ is left-invariant. Hence $\Omega$ can be regarded as a Hermitian-symplectic structure on $(\g,J)$. Using Lemma \ref{noninvHS}, we have that $\Omega$ induces an SKT inner product on
$(\g,J)$. Then Theorem \ref{main2} implies that  ${\mathfrak g}$ is $2$-step nilpotent and therefore ${\mathfrak g}^1$ is contained in the center  $\xi$.  But since $\xi$  is $J$-invariant,   $J (\g^1) \subseteq \xi$ and Theorem \ref{main1} implies the statement.
\end{proof}

Theorem \ref{main3} cannot be directly deduced from Theorem \ref{main1}, since there exist examples of nilpotent Lie algebras not satisfying \eqref{Jxi}.
For instance, we can consider the following

\begin{ex}\label{ex10} {\rm Let $\g$ be  the $10$-dimensional nilpotent Lie algebra with structure equations
$$
\left \{ \begin{array}{l}
d e^j =0, \quad j = 1, \dots, 7,\\
d e^8 = e^1 \wedge e^5 + e^1 \wedge e^6 + e^3 \wedge e^5 + e^3 \wedge e^6,\\
d e^9 = e^2 \wedge e^5 + e^2 \wedge e^6 + e^4 \wedge e^5 + e^4 \wedge e^6,\\
d e^{10} = e^1 \wedge e^8 + e^3 \wedge e^8 + e^2 \wedge e^9 + e^4 \wedge e^9
\end{array}
\right .
$$
(see \cite[Example 3.13]{Rollenske})} endowed with the non-nilpotent complex structure
$$
J e_1 = e_2,  \, J e_3 = e_4, \, J e_5 = e_7, \, J e_8 =  e_9, \, J e_6 = e_{10},
$$
where $\{e_1, \ldots, e_{10}\}$ denotes the dual basis of $\{e^1, \ldots, e^{10}\}$.
We have
$$
\xi = {\mbox {span}} \langle e_7, e_{10}\rangle, \quad {\mathfrak g}^1 = {\mbox {span}} \langle e_8,  e_9, e_{10}\rangle
$$
and thus $J \xi \cap {\mathfrak g}^1 = \{ 0 \}$. In this case a basis of $(1,0)$-forms is then given by
$$
\alpha^1= e^1 + i e^2,  \, \alpha^2= e^3 + i e^4, \, \alpha^3= e^5 + i e^7, \, \alpha^4= e^8 + i e^9, \, \alpha^5 = e^6 + i e^{10}
$$
with complex structure equations
$$
\left \{ \begin{array}{l} d \alpha^j =0, \quad j = 1,2,3,\\
d \alpha^4 = \frac{1}{2} (\alpha^{13}  + \alpha^{1 \overline 3} + \alpha^{23}  + \alpha^{2 \overline 3} +\alpha^{15}  + \alpha^{1 \overline 5} + \alpha^{25}  + \alpha^{2 \overline 5}),\\
d \alpha^5 = \frac{i}{2} (\alpha^{\overline 14}  + \alpha^{1 \overline 4} + \alpha^{\overline 2 4}  + \alpha^{2 \overline 4} ),
\end{array} \right.
$$
where  by $\alpha^{i \overline j}$ we denote $\alpha^i \wedge {\overline \alpha}^j$.
\end{ex}

\section{A classification of   $8$-dimensional  SKT  nilmanifolds}
In this section we classify $8$-dimensional  nilpotent Lie algebras admitting an SKT structure.
We begin describing the   following two families of nilpotent Lie algebras.

\begin{enumerate}
\item[\bf{1.}] \textbf{{\em First family}:} Consider the family of $8$-dimensional  nilpotent  Lie algebras  with  complex structure equations
\begin{equation} \label{structeq1}
\begin{cases}
d \alpha^j= 0, j = 1,2,\\
d \alpha^3 = B_1 \alpha^{1 2} + B_4 \alpha^{1 \overline 1} + B_5 \alpha^{1 \overline 2} + C_3 \alpha^{2 \overline 1} + C_4 \alpha^{2 \overline 2}\,,\\
d \alpha^4 =F_1\alpha^{12} + F_4 \alpha^{1 \overline 1} + F_5 \alpha^{1 \overline 2} + G_3 \alpha^{2 \overline 1} + G_4 \alpha^{2 \overline 2}\,,
\end{cases}
\end{equation}
where   the capitol letters are arbitrary complex numbers and by $\alpha^{i \overline j}$ we denote $\alpha^i \wedge {\overline \alpha}^j$.
Then the inner product  $$g=\sum_{k=1}^4 \a^k\otimes \ov \a^k$$ is SKT if and only if the complex numbers satisfy the equation
\begin{equation} \label{SKTfirstfamily}
\vert B_1 \vert^2 +\vert F_1 \vert^2+  \vert G_3 \vert^2 + \vert B_5 \vert^2 +  \vert C_3 \vert^2+ \vert F_5 \vert^2 =  2 \,  \Re\mathfrak{e}\, (C_4 \overline B_4 + F_4 \overline G_4) \,.
\end{equation}
Notice that in terms of $g$ and of the derivatives of the $\alpha^k$'s, this last equation can be rewritten as
$$
\sum_{j=3}^{4}\left(\|d\alpha^j\|^2+\Re\mathfrak{e}\left[g(d\alpha^j,\alpha^{1\ov 1})g(d\ov \alpha^j,\alpha^{2\ov 2})\right]
-\sum_{k=1}^{2}|g(d\alpha^j,\alpha^{k\ov k})|^2\right)=0
$$
and by a direct computation we have that the Lee form $\theta = J d^* \omega$  is in general  not closed.

Note that if both  $\vert B_4\vert^2 + \vert C_4\vert^2$ and $\vert F_4\vert^2 + \vert G_4\vert^2$ vanish, then the Lie algebra is abelian. So we can suppose for instance that  $\vert B_4\vert^2 + \vert C_4\vert^2 \neq 0$ (otherwise we consider  the change of  basis
$$
\alpha^1 \mapsto \alpha^1, \quad \alpha^2 \mapsto \alpha^2,  \quad \alpha^3 \mapsto \alpha^4,   \quad \alpha^4 \mapsto \alpha^3)
$$
and arguing as in the proof of Lemma 2.14 in \cite{Ugarte} we can change the basis in order to have
$$
\left \{ \begin{array}{l}
d \alpha^j= 0, j = 1,2,\\
d \alpha^3 = \rho \alpha^{1 2} +  \alpha^{1 \overline 1} + B \alpha^{1 \overline 2}  + D \alpha^{2 \overline 2}\,,\\
d \alpha^4 =F_1\alpha^{12} + F_4 \alpha^{1 \overline 1} + F_5 \alpha^{1 \overline 2} + G_3 \alpha^{2 \overline 1} + G_4 \alpha^{2 \overline 2},
\end{array} \right.
$$
with $\rho \in \{ 0,1\}$, $B, D \in \C$.
With the new change
$$
\alpha^1 \mapsto \alpha^1, \quad \alpha^2 \mapsto \alpha^2,  \quad \alpha^3 \mapsto \alpha^3,   \quad \alpha^4 \mapsto \alpha^4 - F_4 \alpha^3,
$$
we then get
$$
\left \{ \begin{array}{l}
d \alpha^j= 0, j = 1,2,\\
d \alpha^3 = \rho \alpha^{1 2} +  \alpha^{1 \overline 1} + B \alpha^{1 \overline 2}  + D \alpha^{2 \overline 2}\,,\\
d \alpha^4 =F'_1\alpha^{12} + F'_5 \alpha^{1 \overline 2} + G'_3 \alpha^{2 \overline 1} + G'_4 \alpha^{2 \overline 2}.
\end{array} \right.
$$
Therefore, if $(F'_1)^2 + (F'_5)^2 + (G'_3)^2 + (G'_4)^2 =0$, the Lie algebra is a direct sum of $\R^2$ with a $6$-dimensional SKT nilpotent Lie algebra. So for instance the direct sum ${\mathfrak h}_3^{\C} \oplus \R^2$  of the $3$-dimensional complex Heisenberg Lie algebra ${\mathfrak h}_3^{\C}$ with $\R^2$ admits an SKT structure. If $(F'_1)^2 + (F'_5)^2 + (G'_3)^2 + (G'_4)^2 \neq 0$,
then $Z_3$ and $Z_4$ belong to the complexification of the center $\xi$, so $\dim \xi  \geq 4$.
Note that $\dim{\mathfrak g}^1 =1$ if and only if  the Lie algebra is isomorphic to the direct sum ${\mathfrak h}_3^{\R}  \oplus \R^5$.
Moreover $\dim {\mathfrak g}^1 =3$ if and only one of the following cases occur
\begin{enumerate}
\item[ a)]  $(F'_1)^2  +  (F'_5)^2 +  (G'_3)^2 \neq 0$, $\rho = B =0$ and $D$  real;

\vspace{0.1cm}
\item[ b)] $F'_1 = F'_5 = G'_3 =0$, $G'_4 \neq 0$ and $\rho^2  + B^2 \neq 0$.
\end{enumerate}

In  case a) we get the Lie algebra
$$
\left \{ \begin{array}{l}
d \alpha^j= 0, j = 1,2,\\
d \alpha^3 = \alpha^{1 \overline 1}   + D \alpha^{2 \overline 2},\\
d \alpha^4 =F'_1\alpha^{12} + F'_5 \alpha^{1 \overline 2} + G'_3 \alpha^{2 \overline 1} + G'_4 \alpha^{2 \overline 2}.
\end{array} \right.
$$
with $(F'_1)^2  +  (F'_5)^2 +  (G'_3)^2 \neq 0$ and $D$  real. In particular, for $D =1$, $ F'_1 = \sqrt{2}$, $F'_5 = G'_3 = G'_4 = 0$ we get the Lie algebra with structure equations
$$
\left \{ \begin{array}{l}
d \alpha^j= 0, j = 1,2,\\
d \alpha^3 = \alpha^{1 \overline 1}   +  \alpha^{2 \overline 2},\\
d \alpha^4 = \sqrt{2} \, \alpha^{12},
\end{array} \right.
$$ which is isomorphic to the direct sum ${\mathfrak h}_7^Q \oplus \R$, where ${\mathfrak h}_7^Q$ is the real $7$-dimensional Lie algebra of Heinseberg type and with $3$-dimensional center.

In  case b)  we get a  Lie algebra  isomorphic to
$$
\left \{ \begin{array}{l}
d \alpha^j= 0,\,\, j = 1,2,\\
d \alpha^3 = \rho \alpha^{1 2} +  \alpha^{1 \overline 1} + B' \alpha^{1 \overline 2},\\
d \alpha^4 =  \alpha^{2 \overline 2},
\end{array} \right.
$$
with $\rho^2 + (B')^2 \neq 0$.

\vspace{0.2cm}
\item[{\bf 2.}] \textbf{{\em Second family}:} Consider the family  of $8$-dimensional  nilpotent Lie algebras equipped with a complex structure having structure equations
\begin{equation} \label{structeq2}
 \begin{cases}
  \begin{aligned}
 d \alpha^j =&0, \quad j = 1,2,3\,,\\
 d \alpha^4 = &F_1 \alpha^{12} + F_2 \alpha^{13}+ F_4 \alpha^{1 \overline 1} + F_5 \alpha^{1 \overline 2} + F_6 \alpha^{1 \overline 3} + G_1 \alpha^{23}   \\[-3pt]
 &+  G_3 \alpha^{2 \overline 1}+ G_4 \alpha^{2 \overline 2}+ G_5 \alpha^{2 \overline 3}
 + H_2 \alpha^{3 \overline 1} + H_3 \alpha^{3 \overline 2} + H_4 \alpha^{3 \overline 3}\,,
\end{aligned}
\end{cases}
\end{equation}
where  the capitol letters are arbitrary complex numbers and  $H_4\neq 0$. In this case requiring that the inner  product $g=\sum_{k=1}^4\alpha^k\otimes \ov \alpha^k$ is  SKT is equivalent to require that the following
equations are satisfied
\begin{equation} \label{secondfamilySKT}
\left \{ \begin{array}{l}
- H_3 \overline F_4  + H_2 \overline G_3 + F_5 \overline F_6 - F_4  \overline G_5 + F_2 \overline F_1=0,\\[4 pt]
-H_3 \overline F_5  + G_4 \overline F_6  + H_2  \overline G_4 - G_3  \overline G_5 + G_1 \overline F_1 =0,\\[4 pt]
-H_4 \overline F_5 + G_5 \overline F_6 + H_2 \overline H_3 -G_3  \overline H_4 + G_1 \overline F_2=0,\\[4 pt]
\vert F_2 \vert^2 + \vert F_6 \vert^2 +\vert H_2 \vert^2= 2 \, {\mbox {Re}} (H_4 \overline F_4),\\[4 pt]
 \vert F_1 \vert^2 + \vert F_5 \vert^2 +  \vert G_3 \vert^2 = 2 \, {\mbox {Re}} (F_4 \overline G_4),\\[4 pt]
  \vert G_1 \vert^2   + \vert G_5 \vert^2 + \vert H_3 \vert^2 = 2 \, {\mbox {Re}} (H_4 \overline G_4 ).
\end{array} \right.
\end{equation}
Every  Lie algebra ${\mathfrak g}$ of this family splits in $\g=V_1\oplus V_2$, where $V_i$ are $J$-invariant  vector subspaces such that
$$\dim V_1=6 ,  \quad  \dim V_2=2,  \quad  [V_1, V_1] \subseteq V_2, \quad V_2\subseteq\xi $$ and there exists a  real $J$-invariant $2$-dimensional vector subspace $V_3$ (generated by $Z_3$ and $\overline Z_3$)  of $V_1$ such that $[V_3,V_3]\neq 0$.  In particular any Lie algebra of this family is $2$-step nilpotent, has $\dim {\mathfrak g}^1 = 2$  and  by a direct computation we have that the Lee form $\theta = J d^* \omega$ is in general  not closed.
\end{enumerate}

Now we are ready to classify $8$-dimensional  nilmanifolds endowed with an invariant complex structure  admitting an SKT metric.

\begin{theorem}\label{classification} Let $M^8 =G/\Gamma $ be an $8$-dimensional  nilmanifold  $($not a torus$)$ with an invariant   complex structure $J$.
There exists an {\rm SKT} metric $g$ on $M$ compatible with $J$ if and only if the Lie algebra $(\g,J)$ belongs to one of the two families described above.
\end{theorem}
\begin{proof} Assume  that there exists an SKT  metric $g$ on $(M^8,J)$. In view of \cite{Ugarte,FG}, we may assume that the metric $g$ is invariant, and so that it is induced by an inner product $g$ on the Lie algebra $\mathfrak g$.

In  the proof of Theorem \ref{main2}
we have already  shown  that there exists  a basis of $(1,0)$-forms $\{ \alpha^1, \alpha^2, \alpha^3, \alpha^4 \}$ such that
$$
\left \{ \begin{array}{l}
d \alpha^j = 0, \quad j = 1, \ldots, p,\\
d \alpha^l \in \Lambda^2 \langle \alpha^1, \ldots, \alpha^p, \overline \alpha^1, \ldots , \overline \alpha^p \rangle, \quad  l = p + 1, \ldots,4,
\end{array} \right.
$$
with $\{ \alpha^1, \ldots, \alpha^p \}$ in the annihilator of $\xi$ and $\{ \alpha^{p + 1}, \ldots, \alpha^4 \}$ in the complexification of  the dual $\xi^*$ of the center.
Applying the  Gram-Schmidt process to the previous basis  $\{  \alpha^1, \dots, \alpha^4 \}$  we get  a unitary coframe
$\{  \tilde \alpha^1, \dots, \tilde \alpha^4 \}$ such that
$$
\left \{ \begin{array}{l}
d \tilde  \alpha^j = 0, \quad j = 1, \ldots, p,\\
d \tilde \alpha^l \in \Lambda^2 \langle \tilde \alpha^1, \ldots, \tilde \alpha^p, \overline {\tilde  \alpha^1}, \ldots , \overline {\tilde \alpha^p} \rangle, \quad  l = p + 1, \ldots,4,
\end{array} \right.
$$
since ${\mbox {span}} \langle  \tilde \alpha^1, \ldots,  \tilde \alpha^j\rangle={\mbox {span}} \langle \alpha^1, \ldots, \alpha^j\rangle $, for any $j= 1, \ldots, 4$. Then it is not restrictive to assume that $ \{\alpha^1,\dots,\alpha^4\}$
is a $g$-unitary coframe, i.e.  that  the fundamental form $\omega$ of $g$ with respect to $\{\alpha^1,\dots,\alpha^4\}$ takes the standard expression:
$$
\omega = -\frac{i} {2} \sum_{j = 1}^4  \alpha^{j}\wedge\a^{ \overline j}.
$$

We can distinguish  three cases according to the dimension $2n - 2p$ of the center $\xi$.

\begin{enumerate}

\item[(a)] If $p =1$, then  $\mathfrak g  \cong {\mathfrak h}_3^{\R} \oplus \R^{5}$ and it belongs to the first family.

\item[(b)] For  $p = 2$  we get  the first family \eqref{structeq1}.

\item[(c)] For  $p = 3$  we obtain  the second  family \eqref{structeq2}.

\end{enumerate}

\end{proof}

\begin{rem} {\rm In dimension $8$  it is  no longer true that the existence of a
strong KT structure on a nilpotent Lie algebra $\mathfrak g$  depends only on the complex structure of  $\mathfrak g$.
Indeed,  for any of two  families of $8$-dimensional  nilpotent Lie algebras with complex structure equations \eqref{structeq1} and  \eqref{structeq2}
consider a generic  $J$-Hermitian metric $g$. The   fundamental form  $\omega$  associated to the Hermitian structure $(J, g)$ can be then expressed as $$
\begin{aligned}
\omega = & a_1 \alpha^{1\ov{1}} + a_2 \alpha^{2\ov{2}} + a_3 \alpha^{3\ov{3}} + a_4 \alpha^{4\ov{4}} +a_5  \alpha^{1 \overline 2}- \overline a_5 \alpha^{2 \overline 1} +a_6 \alpha^{1 \overline 3} -   \overline a_6 \alpha^{3 \overline 1}  \\[-3pt]
& + a_7 \alpha^{1 \overline 4} -\overline  a_7 \alpha^{4 \overline 1}+  a_8 \alpha^{2 \overline 3} - \overline a_{8}\alpha^{3 \overline 2}  +a_9 \alpha^{2 \overline 4} - \overline a_{9} \alpha^{4 \overline 2} +  a_{10}  \alpha^{3\ov{4}} - \ov{a}_{10} \alpha^{4\ov{3}}\,,
\end{aligned}
$$
where  $a_l$, $l = 1, \dots, 10$, are arbitrary   complex numbers     (with  $\overline a_l = - a_l$, for any $l = 1, \ldots, 4$) such that  $\omega$  is positive definite.

For the first family the SKT equation for a generic $J$-Hermitian metric $g$ is:
$$
\begin{aligned}
&- a_3 C_4 \overline B_4 - 2 a_{10} B_4 \overline G_4 - a_3 B_4 \overline C_4 + a_{10} B_1 \overline F_1 + a_3 \vert B_1 \vert^2 - \overline a_{10} \overline B_1 F_1+ a_4 \vert F_5 \vert^2 \\[-3pt]
 &+ a_4 \vert F_1 \vert^2 - \overline a_{10} G_3 \overline C_3 + a_4 \vert G_3 \vert^2 + a_3 \vert B_5 \vert^2 - \overline a_{10} \overline B_5 F_5 - a_4 \overline G_4 F_4 + a_3 \vert C_3 \vert^2\\[-3pt]
  &+ \overline a_{10} \overline C_4 F_4 + a_{10} \overline F_5 B_5 - a_4 \overline F_4 G_4 + \overline a_{10} G_4 \overline B_4
   + a_{10} \overline G_3 C_3 - a_{10} C_4 \overline F_4 =0,
\end{aligned}
$$
so it  is no longer true  that  equation \eqref{SKTfirstfamily} implies that any $J$-Hermitian metric is ${\rm SKT}$.

For the second family the SKT equations  for the generic  $J$-Hermitian metric  $g$ are \eqref{secondfamilySKT}  and so as in the $6$-dimensional case the SKT condition depends only on the complex structure.}
\end{rem}

\smallskip

Nilmanifolds provide also examples of  compact HKT manifolds. We recall that a hyperhermitian  manifold $(M, J_1, J_2, J_3, g)$ is called HKT  if  the fundamental forms $\omega_l (\cdot \, , \cdot \, ) = g( J_l \cdot \, , \cdot \, )$, $l = 1,2,3$, satisfy
$$
J_1 d \omega_1 = J_2 d \omega_2 = J_3 d \omega_3.
$$
This is equivalent to the property that the three Bismut connections associated to the Hermitian
structures $(J_l, g)$ coincide and this connection is said to be an HKT connection. The HKT structure is called strong  or weak   depending on whether the torsion $3$-form of the HKT connection
is closed or not.

By \cite[Theorem 4.2]{BDV}   if  a  $4n$-dimensional nilmanifold  $M^{4n} =  G/\Gamma$ is endowed with
an HKT structure   $(J_1, J_2, J_3,  g)$ induced by a left-invariant HKT structure on $G$,   then the hypercomplex structure has to be abelian, i.e
$$
[J_l X, J_l Y] = [X,Y], \quad l = 1,2,3,
$$
for any $X, Y$ in the Lie algebra ${\mathfrak g}$ of $G$. By \cite{DF} every abelian hypercomplex structure on a (non-abelian) nilpotent Lie algebra gives rise to a weak HKT structure.

In \cite{DF2}   a description of  hypercomplex  $8$-dimensional nilpotent Lie algebras was given. More precisely, it was shown  that  a   $8$-dimensional hypercomplex  nilpotent Lie algebra  ${\mathfrak g}$  has to be $2$-step nilpotent and with   $b_1({\mathfrak g}) \geq 4$.

In particular,  by \cite{DF3}  there exist only three (non-abelian) nilpotent  Lie algebras of dimension $8$  admitting an abelian hypercomplex  structure and they  are abelian extensions of a Lie algebra of Heinsenberg type, i.e. are isomorphic to one of the following:
$$
{\mathfrak h}_5 \oplus \R^3, \quad
{\mathfrak h}_3^{\C} \oplus \R^2, \quad
{\mathfrak h}_7^Q \oplus \R,
$$
where ${\mathfrak h}_5$  is the $5$-dimensional Lie algebra of Heinseberg type with $1$-dimensional center.

As a consequence of Theorem \ref{classification}  we
have the following

\begin{cor} The nilpotent Lie algebra   ${\mathfrak h}_5 \oplus \R^3$ has weak {\rm HKT} structures and it  does not admit any {\rm SKT} structure. The Lie algebras ${\mathfrak h}_3^{\C} \oplus \R^2$ and
${\mathfrak h}_7^Q \oplus \R$ have {\rm SKT} structures and weak {\rm HKT} structures.
\end{cor}

\begin{proof} The nilpotent Lie algebra ${\mathfrak h}_5 \oplus \R^3$ has $1$-dimensional commutator, but  by the description of the two families  a $8$-dimensional SKT nilpotent Lie algebra $\mathfrak g$ with
$\dim {\mathfrak g}^1 = 1$ is isomorphic to ${\mathfrak h}_3^{\R}  \oplus \R^5$.
In the first family there are Lie algebras isomorphic to  ${\mathfrak h}_3^{\C} \oplus \R^2$ and
${\mathfrak h}_7^Q \oplus \R$.
\end{proof}


\begin{thebibliography}{12}


\bibitem{AI} B. Alexandrov, S. Ivanov, Vanishing theorems on Hermitian manifolds, {\em Differ. Geom.
Appl.} {\bf 14} (2001), 251--265.

\bibitem{AG} V. Apostolov and M. Gualtieri, Generalized K\"ahler manifolds with split tangent bundle,
{\em Comm. Math. Phys.} \textbf{271} (2007), 561--575.

\bibitem{BDV}  M. L. Barberis,  I.  Dotti,  M. Verbitsky, Canonical bundles of complex nilmanifolds, with applications to hypercomplex geometry,  {\em Math. Res. Lett.}  {\bf 16} (2009), no. 2, 331--347.

\bibitem{BG} C. Benson and C. S. Gordon, K\" ahler and symplectic structures on nilmanifolds, {\em Topology}
{\bf 27} (1988), 513--518.

\bibitem{Bismut} J.M. Bismut, A local index theorem for non-K\" ahler manifolds, {\em Math. Ann.} {\bf 284} (1989), 681--699.

\bibitem{CG} G. R. Cavalcanti, M. Gualtieri, Generalized complex structures on nilmanifolds,
{\em J. Symplectic Geom.} {\bf 2} (2004), 393--410.

\bibitem{CFGU} L.A. Cordero, M. Fern\' andez, A. Gray, L. Ugarte, Compact nilmanifolds with nilpotent complex structure:  Dolbeault cohomology, {\em Trans. Amer. Math. Soc.}  {\bf 352} (2000), 5405�5433.

\bibitem{donaldson}
S.K. Donaldson, Two-forms on four-manifolds and elliptic equations,
in \emph{Inspired by S.S. Chern}, World Scientific, 2006.


\bibitem{DF3} I. Dotti, A. Fino, Abelian hypercomplex $8$-dimensional nilmanifolds,  {\em Ann. Global Anal. Geom.} {\bf 18}  (2000),  no. 1, 47--59.

\bibitem{DF} I. Dotti, A. Fino, Hyperk\" ahler torsion  structures invariant by nilpotent Lie groups, {\em Classical Quantum Gravity} {\bf 19},  no.3 (2002), 551--562.

\bibitem{DF2} I. Dotti, A. Fino, Hypercomplex $8$-dimensional nilpotent Lie groups,  {\em  J. Pure Applied Algebra} {\bf 184} (2003) no. 1, 41--57.

\bibitem{Enrietti} N. Enrietti, Static SKT metrics on Lie groups,   preprint {\tt arXiv:1009.0620}.

\bibitem{FG} A. Fino, G. Grantcharov, On some properties of the manifolds with skew-symmetric
torsion and holonomy $SU (n)$ and $Sp(n)$,  {\em Advances in Math.} {\bf 189} (2004), 439--450.

\bibitem{FPS} A. Fino, M. Parton, S. Salamon, Families of strong KT structures in six dimensions, {\em Comm. Math. Helv.}  {\bf 79} (2004) no. 2,  317--340.

\bibitem{FT0}   A. Fino, A. Tomassini, Non K\"ahler solvmanifolds with generalized K\"ahler
structure,   {\em J. Symplectic Geom.}  {\bf 7} (2009), no. 2, 1-14.

\bibitem{FT1} A. Fino, A. Tomassini,  Blow ups and resolutions  of strong K\"ahler with torsion metrics,  {\em  Adv. Math.}  {\bf 221}  (2009), no.3, 914--935.

\bibitem{FT2}   A. Fino, A. Tomassini, On astheno-K\"ahler  metrics,  preprint {\tt arXiv:0806.0735}, to appear in Ê{\em J. Lond. Math. Soc}.

\bibitem{GHR}
S. J. Gates, C. M. Hull, M. R\u{o}cek, Twisted multiplets and new supersymmetric nonlinear
sigma models, {\em Nuc. Phys. B} {\bf 248} (1984) 157-186.


\bibitem{Gau}   P. Gauduchon, Hermitian connections and Dirac operators, {\em Boll. Un. Mat. Ital. B}  {\bf 11} (1997), 257--288.

\bibitem{Gau2} P. Gauduchon, La 1-forme de torsione d�une variete hermitienne compacte, {\em Math. Ann. } {\bf 267}  (1984), 495--518.


\bibitem{GK} M. Goze and Y. Khakimdjanov, {\em Nilpotent Lie Algebras}, Mathematics and its Applications {\bf 361},
Kluwer, 1996.

\bibitem{GP} G. Grantcharov, Y.S. Poon, Geometry of hyper-K\"ahler connections with torsion, {\em Comm.
Math. Phys.}  {\bf 213} (2000), 19--37.


\bibitem{Gu} M. Gualtieri, Generalized complex  geometry, DPhil thesis, University of Oxford, 2003,   {\tt arXiv:0401221}.

\bibitem{Has} K. Hasegawa, Minimal models of nilmanifolds, {\em Proc. Amer. Math. Soc.} {\bf 106}  (1989), 65--71.

\bibitem{Ha}
K. Hasegawa, Complex and K\"ahler structures on
compact solvmanifolds, {\em J. Symplectic Geom. } {\bf 3} (2005), 749--767.

\bibitem{Ha2} K. Hasegawa, A note on compact solvmanifolds with K\" ahler
structures, {\em Osaka J. Math.} {\bf 43} (2006), no. 1, 131--135.


\bibitem{Hi2} N.J. Hitchin, Instantons and generalized K\"ahler geometry, {\em Comm. Math. Phys.}
\textbf{265} (2006), 131--164.

\bibitem{HP} P. S. Howe, G. Papadopoulos, Twistor spaces for hyper-K\" ahler manifolds with torsion, {\em Phys.
Lett. B}  {\bf 379}  (1996), 80--86.

\bibitem{LiZhang} T.-J. Li, W. Zhang, Comparing tamed and compatible symplectic cones and
cohomological properties of almost complex manifolds, {\em Comm. Anal. Geom.}, {\bf 17} (2009), no. 4, 651–-683.

\bibitem{MS} T. B. Madsen,  A. Swann, Invariant strong KT geometry on four-dimensional solvable Lie groups, preprint  {\tt arXiv:0911.0535}.

\bibitem{Magnin} L. Magnin, Sur les alg\'ebres de Lie nilpotentes de dimension $\leq 7$, {\em J. Geom. Phys.}  {\bf 3} (1986), 119--144.

\bibitem{Malcev}  A. I. Malcev, On a class of homogeneous spaces, {\em Amer. Math. Soc. Translation Ser. 1} {\bf  9} (1962), 276--307.

\bibitem{Rasmus} R. Mejldal, Complex manifolds and strong geometries with torsion, Master
thesis, Department of Mathematics and Computer Science, University of
Southern Denmark, July 2004.

\bibitem{Mi}  J. Milnor, Curvature of left invariant metrics on Lie groups, {\em Adv. in Math.}  {\bf 21} (1976),  293--329.

\bibitem{Nomizu}  K. Nomizu, On the cohomology of compact homogeneous spaces of nilpotent Lie groups, {\em Ann.  of Math.} {\bf 59} (1954) 531--538.

\bibitem{Rollenske} S. Rollenske, Geometry  of  nilmanifolds  with
left-invariant  complex structure  and
deformations  in the large, {\em Proc. Lond. Math. Soc.} (3) {\bf 99}(2) (2009), 425--460.

\bibitem{RT} F. A. Rossi, A. Tomassini, On strong K\"ahler   and astheno-K\"ahler  metrics  on
nilmanifolds,  preprint (2010).

\bibitem{Salamon} S. Salamon, Complex structures on nilpotent Lie algebras, {\em J. Pure Appl. Algebra }{\bf 157} (2001), 311--333.

\bibitem{StreetsTian} J. Streets, G. Tian, Hermitian Curvature Flow, preprint  {\tt arXiv:0804.4109}, to appear  in {\em JEMS}.


\bibitem{ST} J. Streets,  G. Tian, A Parabolic flow of pluriclosed metrics,  \emph{Int. Math. Res. Notices} (2010) {\bf 2010}, 3101--3133.

\bibitem{strominger} A. Strominger, Superstrings with torsion, {\em Nuclear Phys. B} {\bf 274} (1986) 253--284.

\bibitem{Sw} A. Swann, Twisting Hermitian and hypercomplex geometries,   {\em Duke Math. J.} (2010) {\bf 155},  403--431.

\bibitem{Tomas} A. Tomassini, private communication.

\bibitem{Ugarte} L. Ugarte, Hermitian structures on six dimensional nilmanifolds, {\em Transf. Groups} {\bf 12} (2007),
175--202.

\bibitem{weinkove}
B. Weinkove, The Calabi-Yau equation on almost-K\"ahler four manifolds,
\emph{J. Differential Geom.}  {\bf 76}  (2007),  no. 2, pp. 317--349.

\end{thebibliography}
\end{document}